\documentclass[11pt]{article}

\usepackage{latexsym,color,graphicx,amsmath,amssymb,caption}
\usepackage{enumerate}
\usepackage{lineno}

\usepackage{dsfont}

\def\reals{{\mathbb R}}
\def\eps{{\varepsilon}}

\newtheorem{theorem}{Theorem}[section]

\newtheorem{corollary}[theorem]{Corollary}

\newcommand{\RR}{\ensuremath{\mathbb R}}

\def\eps{{\varepsilon}}

\def\C{{\mathcal C}}

\begin{document}

\def\reals{{\mathbb R}}
\def\cplx{{\mathbb C}}
\def\deg{{\mathsf{deg}}}

\newcommand{\pts}{\mathcal P}
\newcommand{\vrts}{\mathcal V}
\newcommand{\curves}{\mathcal C}
\newcommand{\lines}{\mathcal L}
\newcommand{\motions}{\mathcal M}
\newcommand{\planes}{\Pi}
\newcommand{\flats}{\mathcal F}
\newcommand{\DS}{\ensuremath{\mathbb S}}
\newcommand{\dist}[1]    {{\,{\mathrm{D}_3}\!\left(#1\right)}}
\newcommand{\res}{\text{Res}}

\newcommand{\Deg}{D} 
\makeatletter
\newcommand{\ProofEndBox}{{\ifhmode\unskip\nobreak\hfil\penalty50 \else
          \leavevmode\fi\quad\vadjust{}\nobreak\hfill$\Box$
            \finalhyphendemerits=0 \par}}
\makeatother
\newcommand{\proofend}{\ProofEndBox\smallskip}
\title{Incidences with curves in $\RR^d$\thanks{%
Work on this paper was partially supported by Grant 2012/229 from
the U.S.-Israel Binational Science Foundation, by Grant 892/13 from
the Israel Science Foundation, by the Israeli Centers for Research
Excellence (I-CORE) program (center no.~4/11), and by the Hermann
Minkowski--MINERVA Center for Geometry at Tel Aviv University.
Part of this research was performed while the authors were visiting
the Institute for Pure and Applied Mathematics (IPAM), which is
supported by the National Science Foundation. A preliminary version
of this paper appeared in Proc. European Sympos. Algorithms, 2015,
pages 977--988.} }

\author{
Micha Sharir\thanks{%
Blavatnik School of Computer Science, Tel Aviv University, Tel~Aviv,
Israel; \tt{michas@post.tau.ac.il} } \and
Adam Sheffer\thanks{%
Dept.~of Mathematics, California Institute of Technology, Pasadena,
CA, USA; \tt{adamsh@gmail.com} } \and
Noam Solomon\thanks{%
Blavatnik School of Computer Science, Tel Aviv University, Tel~Aviv,
Israel; \tt{noamsolomon@post.tau.ac.il} } }

\maketitle



\begin{abstract}
We prove that the number of incidences between $m$ points and $n$
bounded-degree curves with $k$ degrees of freedom in $\RR^d$ is
\[
I(m,n) =O\left(m^{\frac{k}{dk-d+1}+\eps}n^{\frac{dk-d}{dk-d+1}}+
\sum_{j=2}^{d-1}
m^{\frac{k}{jk-j+1}+\eps}n^{\frac{d(j-1)(k-1)}{(d-1)(jk-j+1)}}
q_j^{\frac{(d-j)(k-1)}{(d-1)(jk-j+1)}}+m+n\right),
\]
for any $\eps>0$, where the constant of proportionality depends on
$k, \eps$ and $d$, provided that no $j$-dimensional surface of
degree $\le c_j(k,d,\eps)$, a constant parameter depending on $k$,
$d$, $j$, and $\eps$, contains more than $q_j$ input curves, and
that the $q_j$'s satisfy certain mild conditions.

This bound
generalizes a recent result of Sharir and Solomon~\cite{SS14}
concerning point-line incidences in four dimensions (where $d=4$ and
$k=2$), and partly generalizes a recent result of Guth~\cite{Gu14}
(as well as the earlier bound of Guth and Katz~\cite{GK2}) in three
dimensions (Guth's three-dimensional bound has a better dependency
on $q_2$). It also improves a recent $d$-dimensional general
incidence bound by Fox, Pach, Sheffer, Suk, and Zahl~\cite{FPSSZ14},
in the special case of incidences with algebraic curves. Our results
are also related to recent works by Dvir and Gopi~\cite{DG15} and by
Hablicsek and Scherr~\cite{HS14} concerning rich lines in
high-dimensional spaces.
\end{abstract}


\bibliographystyle{plain}


\section{Introduction}

Let $\C$ be a set of curves in $\RR^d$. We say that $\C$
has $k$ \emph{degrees of freedom} with \emph{multiplicity} $s$
if (i) for every $k$ points in $\RR^d$ there are at most $s$
curves of $\cal C$ that are incident to all $k$ points, and (ii)
every pair of curves of $\cal C$ intersect in at most $s$ points.
The bounds that we derive depend more significantly on $k$ than on $s$---see below.

In this paper we derive general upper bounds on the number of
incidences between a set $\pts$ of $m$ points and a set $\C$ of $n$
bounded-degree algebraic curves that have $k$ degrees of freedom
(with some constant multiplicity $s$). We denote the number of these
incidences by $I(\pts,\C)$.

Before stating our results, let us put them in context.
The basic and most studied case involves incidences between points and lines.
In two dimensions, writing $L$ for the given set of $n$ lines, the classical
Szemer\'edi--Trotter theorem~\cite{SzT} yields the worst-case tight bound
\begin{equation} \label{inc2}
I(\pts,L) = O\left(m^{2/3}n^{2/3} + m + n \right) .
\end{equation}
In three dimensions, in the 2010 groundbreaking paper of Guth and Katz~\cite{GK2}, an improved
bound has been derived for $I(\pts,L)$, for a set $\pts$ of $m$ points and a set
$L$ of $n$ lines in $\reals^3$, provided that not too many lines of $L$ lie
in a common plane. Specifically, they showed:
\begin{theorem}[Guth and Katz~\cite{GK2}]
\label {ttt}
Let $\pts$ be a set of $m$ distinct points and $L$ a set of $n$ distinct lines
in $\reals^3$, and let $q_2\le n$ be a parameter,
such that no plane contains more than $q_2$ lines of $L$. Then
$$
I(P,L) = O\left(m^{1/2}n^{3/4} + m^{2/3}n^{1/3}q_2^{1/3} + m + n\right).
$$
\end{theorem}
This bound was a major step in the derivation of the main result of \cite{GK2},
an almost-linear lower bound on the number of distinct distances determined by any
set of $n$ points in the plane, a classical problem posed by Erd{\H o}s in 1946~\cite{Er46}.
Their proof uses several nontrivial tools from algebraic and differential geometry, most
notably the Cayley--Salmon theorem on osculating lines to algebraic surfaces in $\reals^3$,
and additional properties of ruled surfaces. All this machinery comes on top of the main
innovation of Guth and Katz, the introduction of the \emph{polynomial partitioning technique};
see below.

\noindent In four dimensions, Sharir and Solomon~\cite{SS4d} have
obtained the following sharp point-line incidence bound:
\begin{theorem}[Sharir and Solomon~\cite{SS4d}] \label{th:ss4d}
Let $\pts$ be a set of $m$ distinct points and $L$ a set of $n$
distinct lines in $\reals^4$, and let $q_2,q_3\le n$ be parameters,
such that (i) each hyperplane or quadric contains at most $q_3$ lines of $L$,
and (ii) each 2-flat contains at most $q_2$ lines of $L$. Then
\begin {equation}
\label {ma:in0} I(\pts,L) \le 2^{c\sqrt{\log m}} \left( m^{2/5}n^{4/5} + m \right)
+ A\left( m^{1/2}n^{1/2}q_3^{1/4} + m^{2/3}n^{1/3}q_2^{1/3} + n\right) ,
\end {equation}
where $A$ and $c$ are suitable absolute constants.
When $m\le n^{6/7}$ or $m\ge n^{5/3}$, we get the sharper bound
\begin {equation}
\label {ma:in0x} I(\pts,L) \le A \left( m^{2/5}n^{4/5} + m
+ m^{1/2}n^{1/2}q_3^{1/4} + m^{2/3}n^{1/3}q_2^{1/3} + n\right) .
\end {equation}
In general, except for the factor $2^{c\sqrt{\log m}}$, the bound
is tight in the worst case, for any values of $m,n$, with
corresponding suitable ranges of $q_2$ and $q_3$.
\end{theorem}
This improves, in several aspects, an earlier treatment of this problem in
Sharir and Solomon~\cite{SS14}.

Another way to extend the Szemer\'edi--Trotter bound is for curves in the plane with
$k$ degrees of freedom (for lines, $k=2$). This has been done by Pach and Sharir, who showed:\footnote{%
  Their result holds for more general families of curves, not necessarily algebraic, but,
  since algebraicity will be assumed in higher dimensions, we assume it also in the plane.}
\begin{theorem}[Pach and Sharir \cite{PS98}] \label{th:PS}
Let $\pts$ be a set of $m$ points in $\RR^2$ and let $\C$ be a set
of bounded-degree algebraic curves in $\RR^2$ with $k$ degrees of
freedom and with multiplicity $s$. Then
$$
I(\pts,\C) =
O\left(m^{\frac{k}{2k-1}}n^{\frac{2k-2}{2k-1}}+m+n\right) ,
$$
where the constant of proportionality depends on $k$ and $s$.
\end{theorem}
Several special cases of this result, such as the cases of unit
circles and of arbitrary circles, have been considered
separately~\cite{CEGSW,SST}. Unlike the Szemer\'edi-Trotter result
(which arises as a special case of Theorem~\ref{th:PS} with $k=2$),
the bound in Theorem~\ref{th:PS} is not known to be tight for any
$k\ge 3$. In fact, it is known not to be tight for the case of arbitrary circles;
see~\cite{ANP+}.

Here too one can consider the extension of these bounds to higher
dimensions. The literature here is rather scarce, and we only
mention here the work of Sharir, Sheffer and Zahl~\cite{SSZ} on
incidences between points and circles in three dimensions; an
earlier study of this problem by Aronov et al.~\cite{AKS} gives a
different, dimension-independent bound.

The bounds given above include a ``leading term'' that depends only
on $m$ and $n$ (like the term $m^{1/2}n^{3/4}$ in
Theorem~\ref{ttt}), and, except for the two-dimensional case, a
series of ``lower-dimensional'' terms (like the term
$m^{2/3}n^{1/3}q_2^{1/3}$ in Theorem~\ref{ttt} and the terms
$m^{1/2}n^{1/2}q_3^{1/4}$ and $m^{2/3}n^{1/3}q_2^{1/3}$ in
Theorem~\ref{th:ss4d}). The leading terms, in the case of lines,
become smaller as $d$ increases (when $m$ is not too small and not
too large with respect to $n$). Informally, by placing the lines in
a higher-dimensional space, it should become harder to create many
incidences on them.

Nevertheless, this is true only if the setup is ``truly
$d$-dimensional''. This means that not too many lines or curves can
lie in a common lower-dimensional space. The lower-dimensional terms
handle incidences within such lower-dimensional spaces. There is
such a term for every dimension $j=2,\ldots,d-1$, and the
``$j$-dimensional'' term handles incidences within $j$-dimensional
subspaces (which, as the quadrics in the case of lines in four
dimensions in Theorem~\ref{th:ss4d}, are not necessarily linear and
might be algebraic of low constant degree). Comparing the bounds for
lines in two, three, and four dimensions, we see that the
$j$-dimensional term in $d$ dimensions, for $j<d$, is a sharper
variant of the leading term in $j$ dimensions. More concretely, if
that leading term in $j$ dimensions is $m^an^b$ then its counterpart
in the $d$-dimensional bound is of the form $m^an^tq_j^{b-t}$, where
$q_j$ is the maximum number of lines that can lie in a common
$j$-dimensional flat or low-degree variety,
and $t$ depends on $j$ and $d$.

\noindent{\bf Our results.} In this paper we consider a grand
generalization of these results, to the case where $\C$ is a family
of bounded-degree algebraic curves with $k$ degrees of freedom (and
some multiplicity $s$) in $\RR^d$. This is a very ambitious and
difficult project, and the challenges that it faces seem to be
enormous. Here we make the first, and fairly significant, step in
this direction, and obtain the following bounds. As the exponents in
the bounds are rather cumbersome expressions in $d$, $k$, and $j$,
we first state the special case of $d=3$ (and prove it separately),
and then give the general bound in $d$ dimensions.

\begin{theorem}[Curves in $\RR^3$] \label{th:mainInc3}
Let $k\ge 2$ be an integer, and let $\eps>0$. Then
there exists a constant $c(k,\eps)$ that depends on $k$ and $\eps$,
such that the following holds.
Let $\pts$ be a set of $m$ points and $\C$ a set of $n$ irreducible
algebraic curves of constant degree with $k$ degrees of freedom
(and some multiplicity $s$) in $\RR^3$, such that every algebraic
surface of degree at most $c(k,\eps)$ contains at most $q_2$ curves
of $\C$. Then
\[
I(\pts,\C) =O\left(m^{\frac{k}{3k-2}+\eps}n^{\frac{3k-3}{3k-2}}+
m^{\frac{k}{2k-1}+\eps}n^{\frac{3k-3}{4k-2}}q_2^{\frac{k-1}{4k-2}}+m+n\right),
\]
where the constant of proportionality depends on $k$, $s$, and $\eps$
(and on the degree of the curves).
\end{theorem}
The corresponding result in $d$ dimensions is as follows.
\begin{theorem}[Curves in $\RR^d$] \label{th:IncD}
Let $d\ge 3$ and $k\ge 2$ be integers, and let $\eps>0$. Then there
exist constants $c_j(k,d,\eps)$, for $j=2,\ldots,d-1$, that depend
on $k$, $d$, $j$, and $\eps$, such that the following holds. Let
$\pts$ be a set of $m$ points and $\C$ a set of $n$ irreducible
algebraic curves of constant degree with $k$ degrees of freedom (and
some multiplicity $s$) in $\RR^d$. Moreover, assume that, for
$j=2,\ldots, d-1$, every $j$-dimensional algebraic variety of degree
at most $c_j(k,d,\eps)$ contains at most $q_j$ curves of $\C$, for
given parameters $q_2 \le \ldots \le q_{d-1}\le n$. Then we have
\[
I(\pts,\C) =O\left(m^{\frac{k}{dk-d+1}+\eps}n^{\frac{dk-d}{dk-d+1}}+
\sum_{j=2}^{d-1}
m^{\frac{k}{jk-j+1}+\eps}n^{\frac{d(j-1)(k-1)}{(d-1)(jk-j+1)}}
q_j^{\frac{(d-j)(k-1)}{(d-1)(jk-j+1)}}+m+n\right),
\]
where the constant of proportionality depends on $k$, $s$, $d$, and $\eps$
(and on the degree of the curves),
provided that, for any $2\le j<l \le d$, we have (with the
convention that $q_d=n$)
\begin {equation} \label {eq:ontheq}
q_j \ge \left(\frac {q_{l-1}} {q_l}\right)^{l(l-2)}q_{l-1} .
\end {equation}
\end{theorem}

%

\noindent{\bf Discussion.} The advantages of our results are
obvious: They provide the first nontrivial bounds for the general
case of curves with any number of degrees of freedom in any
dimension (with the exception of one previous study of Fox et
al.~\cite{FPSSZ14}, in which weaker bounds are obtained, for
arbitrary varieties instead of algebraic curves). Apart for the
$\eps$ in the exponents, the leading term is ``best possible,'' in
the sense that (i) the polynomial partitioning technique~\cite{GK2}
that our analysis employs (and that has been used in essentially all recent
works on incidences in higher dimensions) yields a recurrence that
solves to this bound, and, moreover,
(ii) it is (nearly) worst-case tight for \emph{lines} in two,
three, and four dimensions (as shown in the respective works cited
above), and in fact is likely to be tight for lines in higher
dimensions too, using a suitable extension of a construction, due to Elekes
and used in \cite{GK2,SS4d}.

Nevertheless, our bounds are not perfect, and tightening them further is a major challenge for
future research. Specifically:

\noindent{\bf (i)} The bounds involve the factor $m^\eps$. As the existing
works indicate, getting rid of this factor is no small feat. Although the factor does not show up in the
cases of lines in two and three dimensions, it already shows up (sort of) in four dimensions
(Theorem~\ref{th:ss4d}), as well as in the case of circles in three dimensions~\cite{SSZ}.
(A recent study of Guth~\cite{Gu14} also pays this factor for the case of lines in three
dimensions, in order to simplify the original analysis in the Guth--Katz paper~\cite{GK2}.)
See the proofs and comments below for further elaboration of this issue.

\noindent{\bf (ii)} The condition that no surface of degree $c_j(k,d,\eps)$ contains too many curves
of $\C$, for $j=2,\ldots,d-1$, is very restrictive, especially since the actual values
of these constants that arise in the proofs can be quite large. Again, earlier works
also ``suffer'' from this handicap, such as Guth's work~\cite{Gu14} mentioned above,
as well as an earlier version of Sharir and Solomon's four-dimensional bound~\cite{SS14}.

\noindent{\bf (iii)} Finally, the lower-dimensional terms that we obtain are not best possible. For example,
the bound that we get in Theorem~\ref{th:mainInc3} for the case of lines in $\RR^3$ ($k=2$) is
$O(m^{1/2+\eps}n^{3/4} + m^{2/3+\eps}n^{1/2}q_2^{1/6}+m+n)$. When $q_2 \ll n$, the
two-dimensional term $m^{2/3+\eps}n^{1/2}q_2^{1/6}$ in that bound is worse than the corresponding term
$m^{2/3}n^{1/3}q_2^{1/3}$ in Theorem~\ref{ttt} (even when ignoring the factor $m^\eps$).

Our results are also related to recent works by Dvir and Gopi~\cite{DG15} and by Hablicsek and Scherr~\cite{HS14},
that study rich lines in high dimensions. Specifically, let $\pts$ be a set of $n$ points in $\reals^d$
and let $L$ be a set of $r$-rich lines (each line of $L$ contains at least $r$ points of $\pts$).
If $|L|=\Omega(n^2/r^{d+1})$ then there exists a hyperplane containing $\Omega(n/r^{d-1})$ points of $\pts$.
Our bounds might be relevant for extending this result to rich curves. Concretely,
for a set $\pts$ of $n$ points in $\reals^d$ and a collection $\C$ of $r$-rich constant-degree
algebraic curves, if $|\C|$ is too large
then the incidence bound becomes larger than our ``leading term'', indicating that some lower-dimensional
surface must contain many curves of $\C$, from which it might be possible to also deduce that such a surface
has to contain many points of $\pts$. While such an extension is not straightforward, we believe that
it is doable, and plan to investigate it in our future work.

As in the classical work of Guth and Katz~\cite{GK2}, and in the
numerous follow-up studies of related problems, here too we use the
polynomial partitioning method, as pioneered in \cite{GK2}. The
reason why our bounds suffer from the aforementioned handicaps is
that we use a partitioning polynomial of (large but) constant
degree. (The idea of using constant-degree partitioning polynomials for problems of this kind
is due to Solymosi and Tao~\cite{ST11}.) When using a polynomial of a larger, non-constant degree, we
face the difficult task of bounding incidences between points and
curves that are fully contained in the zero set of the polynomial,
where the number of curves of this kind can be large, because the
polynomial partitioning technique has no control over this value. We
remark that for lines we have the classical Cayley--Salmon theorem
(see, e.g., Guth and Katz~\cite{GK2}), which essentially bounds the
number of lines that can be fully contained in an algebraic surface
of a given degree, unless the surface is ruled by lines. However,
such a property has not been known for more general curves.
Nevertheless, Nilov and Skopenkov~\cite{NS12} have recently
established such a result involving lines and circles
in $\reals^3$, and, very recently, Guth and Zahl~\cite{GZ}
have done the same for general algebraic curves in three dimensions.
Handling these incidences
requires heavy-duty machinery from algebraic geometry, and leads to
profound new problems in that domain that need to be tackled.

In contrast, using a polynomial of constant degree makes this part
of the analysis much simpler, as can be seen below, but then
handling incidences within the cells of the partition becomes
non-trivial, and a naive approach yields a bound that is too large.
To handle this part, one uses induction within each cell of the
partitioning, and it is this induction process that is responsible
for the weaker aspects of the (lower-dimensional terms in the)
resulting bound. Nevertheless, with these ``sacrifices'' we are able
to obtain a ``general purpose'' bound that holds for a broad
spectrum of instances. It is our hope that this study will motivate
further research on this problem that would improve our results
along the ``handicaps'' mentioned earlier.
Recalling how inaccessible were
these kinds of problems prior to Guth and Katz's breakthroughs seven and five
years ago, it is quite gratifying that so much new ground can be
gained in this area, including the progress made in this paper.


\noindent{\bf Background.} Incidence problems have been a major
topic in combinatorial and computational geometry for the past
thirty years, starting with the aforementioned Szemer\'edi-Trotter
bound \cite{SzT} back in 1983 (and even earlier). Several
techniques, interesting in their own right, have been developed, or
adapted, for the analysis of incidences, including the
crossing-lemma technique of Sz\'ekely~\cite{Sz}, and the use of
cuttings as a divide-and-conquer mechanism (e.g., see~\cite{CEGSW}).
Connections with range searching and related algorithmic problems in
computational geometry have also been noted and exploited, and
studies of the Kakeya problem (see, e.g., \cite{T}) indicate the
connection between this problem and incidence problems. See Pach and
Sharir~\cite{PS} for a comprehensive (albeit a bit outdated) survey
of the topic.

The landscape of incidence geometry has dramatically changed in the
past seven years, due to the infusion, in two groundbreaking papers
by Guth and Katz~\cite{GK,GK2}, of new tools and techniques drawn
from algebraic geometry. Although their two direct goals have been to obtain
a tight upper bound on the number of joints in a set of lines in three
dimensions \cite{GK}, and a near-linear lower bound for the classical distinct
distances problem of Erd{\H o}s \cite{GK2}, the new tools have quickly
been recognized as useful for incidence bounds.
See \cite{EKS,KMSS,KMS,SSZ,ST11,Za1,Za2} for a sample of recent
works on incidence problems that use the new algebraic machinery.

The present paper continues this line of research, and aims at extending
the collection of instances where nontrivial incidence bounds in higher
dimensions can be obtained.


\section{The three-dimensional case}

{\bf Proof of Theorem \ref{th:mainInc3}.} We fix $\eps > 0$, and
prove by induction on $m+n$ that
\begin{equation} \label{eqind}
I(\pts,\C) \le \alpha_{1} \left( m^{\frac{k}{3k-2}+\eps}n^{\frac{3k-3}{3k-2}}+m^{\frac{k}{2k-1}+\eps}n^{\frac{3k-3}{4k-2}}q^{\frac{k-1}{4k-2}}
\right) +\alpha_{2}(m+n),
\end{equation}
where $\alpha_{1},\alpha_{2}$ are sufficiently large constants,
$\alpha_1$ depends on $\eps$ and $k$ (and $s$), and $\alpha_2$
depends on $k$ (and $s$).

For the induction basis, the case where $m,n$ are sufficiently small
constants can be handled by choosing sufficiently large values of
$\alpha_{1},\alpha_{2}$.

Since the incidence graph, as a subgraph of $\pts\times\C$, does not
contain $K_{k,s+1}$ as a subgraph, the K\H ov\' ari-S\'os-Tur\'an
theorem (e.g., see \cite[Section 4.5]{Mat02}) implies that
$I(\pts,\C) = O(mn^{1-1/k} + n)$, where the constant of proportionality
depends on $k$ (and $s$). When $m=O(n^{1/k})$, this implies the
bound $I(\pts,\C) = O(n)$, which is subsumed in (\ref{eqind}) if we
choose $\alpha_2$ sufficiently large. We may thus assume that
$n\le cm^k$, for some absolute constant $c$.

\noindent {\bf Applying the polynomial partitioning technique.} We
construct an \emph{$r$-partitioning polynomial} $f$ for $\pts$, for
a sufficiently large constant $r$ (depending on $\eps$). That is, as
established in Guth and Katz~\cite{GK2}, $f$ is of degree
$O(r^{1/3})$ (the constant in the $O$ notation is an absolute
constant), and the complement of its zero set $Z(f)$ is partitioned
into $u=O(r)$ open connected cells, each containing at most $m/r$
points of $\pts$.
%
Denote the (open) cells of the partition as $\tau_1, \ldots, \tau_u$.
For each $i=1,\ldots,u$, let $\C_i$
denote the set of curves of $\C$ that intersect $\tau_i$ and let
$\pts_i$ denote the set of points that are contained in $\tau_i$. We
set $m_i=|\pts_i|$ and $n_i=|\C_i|$, for $i=1,\ldots,u$, and $m' = \sum_i m_i$, and
notice that $m_i\le m/r$ for each $i$ (and $m'\le m$). An obvious property
(which is a consequence of B\'ezout's theorem, see, e.g., \cite[Theorem A.2]{ST11})
is that every curve of $\C$ intersects $O(r^{1/3})$ cells of $\RR^3\setminus Z(f)$.
Therefore, $\sum_i n_i \le bnr^{1/3}$, for a suitable absolute constant $b>1$
(that depends on the degree of the curves in $\C$). Using H\"older's inequality, we have
\begin{align*}
\sum_i n_i^{\frac{3k-3}{3k-2}} &\le \left(\sum_i n_i\right)^{\frac{3k-3}{3k-2}}
\left(\sum_i 1\right)^{\frac{1}{3k-2}} \le b'\left(nr^{\frac{1}{3}}\right)^{\frac{3k-3}{3k-2}}r^{\frac{1}{3k-2}}
= b'n^{\frac{3k-3}{3k-2}}r^{\frac{k}{3k-2}}, \\[2mm]
\sum_i n_i^{\frac{3k-3}{4k-2}} &\le \left(\sum_i
n_i\right)^{\frac{3k-3}{4k-2}} \left(\sum_i 1\right)^{\frac{k+1}{4k-2}} \le
b'\left(nr^{\frac{1}{3}}\right)^{\frac{3k-3}{4k-2}}r^{\frac{k+1}{4k-2}}
= b'n^{\frac{3k-3}{4k-2}}r^{\frac{k}{2k-1}},
\end{align*}
for another absolute constant $b'$.
Combining the above with the induction hypothesis, applied
within each cell of the partition, implies
$$
\sum_i I ( \pts_i,\C_i) \le \sum_i \left(\alpha_{1} \left( m_i^{\frac{k}{3k-2}+\eps}n_i^{\frac{3k-3}{3k-2}}+
m_i^{\frac{k}{2k-1}+\eps}n_i^{\frac{3k-3}{4k-2}}q_2^{\frac{k-1}{4k-2}}\right)+\alpha_{2}(m_i+n_i)\right)
$$
$$
\le \alpha_{1}\left(\frac{m^{\frac{k}{3k-2}+\eps}}{r^{\frac{k}{3k-2}+\eps}}
\sum_i n_i^{\frac{3k-3}{3k-2}}
+ \frac{m^{\frac{k}{2k-1}+\eps}q_2^{\frac{k-1}{4k-2}}}{r^{\frac{k}{2k-1}+\eps}}
\sum_i n_i^{\frac{3k-3}{4k-2}}\right) + \sum_i\alpha_{2}(m_i+n_i)
$$
$$
\le \alpha_{1}b'\left(
\frac{m^{\frac{k}{3k-2}+\eps}n^{\frac{3k-3}{3k-2}}}{r^{\eps}} +
\frac{m^{\frac{k}{2k-1}+\eps}n^{\frac{3k-3}{4k-2}}q_2^{\frac{k-1}{4k-2}}}{r^{\eps}}\right)
+ \alpha_{2}\left(m'+bnr^{1/3}\right).
$$
Our assumption that $n=O(m^k)$ implies that
$n=O\left(m^{\frac{k}{3k-2}}n^{\frac{3k-3}{3k-2}}\right)$ (with an
absolute constant of proportionality). Thus, when $\alpha_{1}$ is
sufficiently large with respect to $r,k$, and $\alpha_{2}$, we have
\[ \sum_i I(\pts_i,\C_i) \le 2\alpha_{1}b'\left(\frac{m^{\frac{k}{3k-2}+\eps}n^{\frac{3k-3}{3k-2}}}{r^{\eps}}   +
\frac{m^{\frac{k}{2k-1}+\eps}n^{\frac{3k-3}{4k-2}}q_2^{\frac{k-1}{4k-2}}}{r^{\eps}} \right) + \alpha_{2}m'. \]

When $r$ is sufficiently large, such that $r^\eps \ge 6b'$, we have
\begin{equation}\label{eq:incCells}
\sum_i I(\pts_i,\C_i) \le \frac{\alpha_{1}}{3} \left(m^{\frac{k}{3k-2}+\eps}n^{\frac{3k-3}{3k-2}} +
m^{\frac{k}{2k-1}+\eps}n^{\frac{3k-3}{4k-2}}q_2^{\frac{k-1}{4k-2}} \right)+  \alpha_{2}m'.
\end{equation}

\noindent {\bf Incidences on the zero set $Z(f)$.} It remains to
bound incidences with points that lie on $Z(f)$. Set $\pts_0:= \pts \cap Z(f)$
and $m_0=|\pts_0|=m-m'$. Let $\C_0$ denote the set of
curves that are fully contained in $Z(f)$, and set $\C': = \C \setminus \C_0$,
$n_0 := |\C_0|$, and $n' := |\C'| =n-n_0$. Since every curve of $\C'$ intersects $Z(f)$ in
$O(r^{1/3})$ points, we have, taking $\alpha_1$ to be sufficiently large, and arguing as above,
\begin{equation} \label{eq:inc'}
I(\pts_0,\C') = O(nr^{1/3}) \le \frac{\alpha_1}{3}
m^{\frac{k}{3k-2}+\eps}n^{\frac{3k-3}{3k-2}} .
\end{equation}
Finally, we consider the number of incidences between points of
$\pts_0$ and curves of $\C_0$. For this, we set $c(k,\eps)$ to be
the degree of $f$, which is $O(r^{1/3})$, and can be taken to be $O((6b')^{1/(3\eps)})$.
Then, by the assumption of the theorem, we have $|\C_0| \le q_2$. We
consider a generic plane $\pi \subset \RR^3$ and project $\pts_0$
and $\C_0$ onto two respective sets $\pts^*$ and $\C^*$ on $\pi$.
Since $\pi$ is chosen generically, we may assume that no two points
of $\pts_0$ project to the same point in $\pi$, and that no pair of
distinct curves in $\C_0$ have overlapping projections in $\pi$.
Moreover, the projected curves still have $k$ degrees of freedom, in
the sense that, given any $k$ points on the projection $\gamma^*$ of
a curve $\gamma\in\C_0$, there are at most $s-1$ other projected
curves that go through all these points. This is argued by lifting
each point $p$ back to the point $\bar{p}$ on $\gamma$ in
$\reals^3$, and by exploiting the facts that the original curves
have $k$ degrees of freedom, and that, for a sufficiently generic
projection, any curve that does not pass through $\bar{p}$ does not contain any
point that projects to $p$. The number of intersection points
between a pair of projected curves may increase
but it must remain a
constant since these are intersection points between constant-degree
algebraic curves with no common components. By applying Theorem
\ref{th:PS}, we obtain
\begin{equation*}
I(\pts_0,\C_0) = I(\pts^*,\C^*) = O(m_0^{\frac{k}{2k-1}}q_2^{\frac{2k-2}{2k-1}} + m_0 + q_2),
\end{equation*}
where the constant of proportionality depends on $k$ (and $s$).
Since $q_2\le n$ and $m_0\le m$, we have
${\displaystyle m_0^{\frac{k}{2k-1}}q_2^{\frac{2k-2}{2k-1}} \le
m^{\frac{k}{2k-1}}n^{\frac{3k-3}{4k-2}}q_2^{\frac{k-1}{4k-2}}}$. We thus get
that $I(\pts_0,\C_0)$ is at most
\begin{equation} \label{eq:inc0}
O\left(m^{\frac{k}{2k-1}}n^{\frac{3k-3}{4k-2}}q_2^{\frac{k-1}{4k-2}} +
n + m_0\right)
\le \frac{\alpha_1}{3}m^{\frac{k}{2k-1}}n^{\frac{3k-3}{4k-2}}q_2^{\frac{k-1}{4k-2}} +
b_2n + \alpha_2 m_0 ,
\end{equation}
for sufficiently large $\alpha_1$ and $\alpha_2$; the constant $b_2$ comes from
Theorem \ref{th:PS}, and is independent of $\eps$ and of the choices for $\alpha_1,\alpha_2$ made so far.

By combining \eqref{eq:incCells}, \eqref{eq:inc'}, and \eqref{eq:inc0}, including the case $m=O(n^{1/k})$, and choosing $\alpha_2$ sufficiently large, we obtain
\begin{equation*}
I(\pts,\C) \le \alpha_{1} \left( m^{\frac{k}{3k-2}+\eps}n^{\frac{3k-3}{3k-2}}+m^{\frac{k}{2k-1}+\eps}n^{\frac{3k-3}{4k-2}}q_2^{\frac{k-1}{4k-2}}
\right) +\alpha_{2}(m+n).
\end{equation*}
This completes the induction step and thus the proof of the theorem.
\proofend

\noindent{\bf Example 1: The case of lines.} Lines in $\reals^3$
have $k=2$ degrees of freedom, and we almost get the bound of Guth
and Katz in Theorem~\ref{ttt}. There are three differences that make this derivation somewhat inferior
to that in Guth and Katz~\cite{GK2}, as detailed in items (i)--(iii) in the discussion in the introduction.
%
We also recall the two follow-up studies of point-line incidences in $\RR^3$, of
Guth~\cite{Gu14} and of Sharir and Solomon~\cite{SS3d}. Guth's bound
suffers from weaknesses (i) and (ii), but avoids (iii), using a
fairly sophisticated inductive argument. Sharir and Solomon's bound
avoids (i) and (iii), and almost avoids (ii), in a sense that we do
not make explicit here. In both cases, considerably more
sophisticated machinery is needed to achieve these improvements.

\noindent{\bf Example 2: The case of circles.}
Circles in $\reals^3$ have $k=3$ degrees of freedom, and we get the bound
\[ I(\pts,\C) =O\left(m^{3/7+\eps}n^{6/7}+
m^{3/5+\eps}n^{3/5}q_2^{1/5}+m+n\right).\]
The leading term is the same as in Sharir et al.~\cite{SSZ}, but the second term is weaker,
because it relies on the general bound of Pach and Sharir (Theorem~\ref{th:PS}),
whereas the bound in \cite{SSZ} exploits an
improved bound for point-circle incidences, due to
Aronov et al.~\cite{AKS}, which holds in any dimension. If we plug that bound into the above
scheme, we obtain an exact reconstruction of the bound in \cite{SSZ}.
In addition, considering the items (i)--(iii) discussed earlier, we note:
(i) The requirements in \cite{SSZ} about the maximum number of circles on a surface
are weaker, and are only for planes and spheres.
(ii) The $m^\eps$ factors are present in both bounds. (iii) Even after the improvement
noted above, the bounds still seem to be weak in terms of their dependence on $q_2$, and
improving this aspect, both here and in \cite{SSZ}, is a challenging open problem.


Theorem \ref{th:mainInc3} can easily be restated as bounding the
number of \emph{rich points}.
\begin{corollary} \label{co:rich}
For each $\eps>0$ there exists a parameter $c(k,\eps)$ that depends
on $k$ and $\eps$, such that the following holds. Let $\C$ be a set of $n$
irreducible algebraic curves of constant degree and with $k$ degrees
of freedom (with some multiplicity $s$) in $\RR^3$. Moreover, assume
that every surface of degree at most $c(k,\eps)$ contains at most
$q_2$ curves of  $\C$. Then, there exists some constant
$r_0(k,\eps)$ depending on $\eps, k$ (and $s$), such that for any
$r\ge r_0(k,\eps)$, the number of points that are incident to at
least $r$ curves of $\C$ (so-called \emph{$r$-rich points}), is\\
${\displaystyle O\left( \frac{n^{3/2+\eps}}{r^{\frac{3k-2}{2k-2}+\eps}} +
\frac{n^{3/2+\eps}q_2^{1/2+\eps}}{r^{\frac{2k-1}{k-1}+\eps}} + \frac{n}{r} \right)}$,
where the constant of proportionality depends on $k$, $s$ and
$\eps$.
\end{corollary}

\noindent{\bf Proof.}
Denoting by $m_r$ the number of $r$-rich points, the corollary is obtained by
combining the upper bound in Theorem~\ref{th:mainInc3} with the lower bound $rm_r$.
\proofend


\section{Incidences in higher dimensions}

\noindent{\bf Proof of Theorem~\ref{th:IncD}.} Again, we fix
$\eps>0$, and prove, by double induction, where the outer induction
is on the dimension $d$ and the inner induction is on $m+n$, that
$I(\pts,\C)$ is at most
\begin{equation}
\label{eq:indgend}
\alpha_{1,d} \left(m^{\frac{k}{dk-d+1}+\eps}n^{\frac{dk-d}{dk-d+1}}+
\sum_{j=2}^{d-1}
m^{\frac{k}{jk-j+1}+\eps}n^{\frac{d(j-1)(k-1)}{(d-1)(jk-j+1)}}q_j^{\frac{(d-j)(k-1)}{(d-1)(jk-j+1)}}\right)
+\alpha_{2,d}(m+n),
\end{equation}
where $\alpha_{1,d},\alpha_{2,d}$ are sufficiently large constants,
$\alpha_{1,d}$ depends on $k, \eps, d$ (and $s$), and $\alpha_{2,d}$
depends only on $d, k$ (and $s$).

For the outer induction basis, the case $d=2$ follows by
Theorem~\ref{th:PS}, and the case $d=3$ is just
Theorem~\ref{th:mainInc3}, proved in the previous section. We assume
therefore that the claim holds up to dimension $d-1$, and prove it
in dimension $d\ge 4$. The base case of the inner induction (that
is, when $d$ is fixed, we induct over $m+n$) is when $m,n$ are
sufficiently small constants. The bound in (\ref{eq:indgend}) can
then be enforced by choosing sufficiently large values of
$\alpha_{1,d},\alpha_{2,d}$.

The case $m=O(n^{1/k})$ can be handled exactly as for $d=3$, so we may
assume, as before, that $n\le c m^k$ for some absolute constant $c$.

\noindent {\bf Applying the polynomial partitioning technique.}
The analysis is somewhat repetitive and resembles the one in the
previous section, although many details are different; it is given
in detail for the convenience of the reader.

Let $f$ be an $r$-partitioning polynomial, for a sufficiently large constant $r$.
According to the polynomial partitioning theorem \cite{GK2}, we have
$\deg f = O(r^{1/d})$. Denote the (open) cells of the partition as
$\tau_1, \ldots, \tau_u$, where
$u=O(r)$. For each $i=1,\ldots,u$, let $\C_i$
denote the set of curves of $\C$ that intersect $\tau_i$ and let
$\pts_i$ denote the set of points that are contained in $\tau_i$. We
set $m_i=|\pts_i|$, and $n_i=|\C_i|$, for $i=1,\ldots,u,$ and $m' =
\sum_i m_i$, and notice that $m_i\le m/r$ for each $i$ (and $m' \le
m$). Arguing as before, every curve of $\C$ intersects at most
$\deg(f) = O(r^{1/d})$ cells of $\RR^d\setminus Z(f)$. Therefore,
$\sum_i n_i \le b_dnr^{1/d}$, for a suitable constant $b_d>1$ that
depends on $d$ and the degree of the curves.
Using H\"older's inequality, we have
\begin{align*}
\sum_i n_i^{\frac{dk-d}{dk-d+1}} &\le
b'_d\left(nr^{\frac{1}{d}}\right)^{\frac{dk-d}{dk-d+1}}r^{\frac{1}{dk-d+1}}
\le b'_dn^{\frac{dk-d}{dk-d+1}}r^{\frac{k}{dk-d+1}}, \quad \text{and} \\
\sum_i n_i^{\frac{d(j-1)(k-1)}{(d-1)(jk-j+1)}} &\le
b'_d\left(nr^{\frac{1}{d}}\right)^{\frac{d(j-1)(k-1)}{(d-1)(jk-j+1)}}r^{\frac{dk-jk+j-1}{(d-1)(jk-j+1)}}
\le b'_dn^{\frac{d(j-1)(k-1)}{(d-1)(jk-j+1)}}r^{\frac{k}{jk-j+1}} ,
\end{align*}
for each $j=2,\ldots,d-1$, where $b'_d$ is another constant
parameter that depends on $d$. Combining the above with the
induction hypothesis implies that $\sum_i I ( \pts_i,\C_i)$ is at
most
$$
\sum_i \left(\alpha_{1,d} \left(
m_i^{\frac{k}{dk-d+1}+\eps}n_i^{\frac{dk-d}{dk-d+1}}+\sum_{j=2}^{d-1}
m_i^{\frac{k}{jk-j+1}+\eps} n_i^{\frac{d(j-1)(k-1)}{(d-1)(jk-j+1)}}q_j^{\frac{(d-j)(k-1)}{(d-1)(jk-j+1)}} \right)
+ \alpha_{2,d}(m_i+n_i)\right)
$$
$$
\le \alpha_{1,d}\left(\frac{m^{\frac{k}{dk-d+1}+\eps}}{r^{\frac{k}{dk-d+1}+\eps}}
\sum_i n_i^{\frac{dk-d}{dk-d+1}} + \sum_{j=2}^{d-1}
\frac{m^{\frac{k}{jk-j+1}+\eps}q_j^{\frac{(d-j)(k-1)}{(d-1)(jk-j+1)}}}{r^{\frac{k}{jk-j+1}+\eps}}
\sum_i n_i^{\frac{d(j-1)(k-1)}{(d-1)(jk-j+1)}}\right) + \sum_i\alpha_{2,d}(m_i+n_i)
$$
$$
\le \alpha_{1,d} b'_d\left(
\frac{m^{\frac{k}{dk-d+1}+\eps}n^{\frac{dk-d}{dk-d+1}}}{r^{\eps}} +
\frac{\sum_{j=2}^{d-1} m^{\frac{k}{jk-j+1}+\eps} n^{\frac{d(j-1)(k-1)}{(d-1)(jk-j+1)}}
q_j^{\frac{(d-j)(k-1)}{(d-1)(jk-j+1)}}}{r^\eps} \right)
+ \alpha_{2,d}\left(m'+b_dnr^{1/d}\right).
$$
Since we assume that $n=O(m^k)$, we have
$n=O\left(m^{\frac{k}{dk-d+1}}n^{\frac{dk-d}{dk-d+1}}\right)$,
with a constant of proportionality that depends only on $d$.
Thus, when $\alpha_{1,d}$ is sufficiently large with
respect to $r,d$, and $\alpha_{2,d}$, we have
\[
\sum_i I(\pts_i,\C_i) \le 2\alpha_{1,d}
b\left(\frac{m^{\frac{k}{dk-d+1}+\eps}n^{\frac{dk-d}{dk-d+1}}}{r^{\eps}}
+
\frac{m^{\frac{k}{2k-1}+\eps}n^{\frac{dk-d}{(d-1)(2k-1)}}q^{\frac{(k-1)(d-2)}{(d-1)(2k-1)}}}{r^{\eps}}
\right) + \alpha_{2,d}m'.
\]
When $r$ is sufficiently large, such that $r^\eps \ge 6b'$, the bound is at most
\begin{equation}\label{eq:incCellsd}
\frac{\alpha_{1,d}}{3} \left(m^{\frac{k}{dk-d+1}+\eps}n^{\frac{dk-d}{dk-d+1}} +
\sum_{j=2}^{d-1} m^{\frac{k}{jk-j+1}+\eps} n^{\frac{d(j-1)(k-1)}{(d-1)(jk-j+1)}}
q_j^{\frac{(d-j)(k-1)}{(d-1)(jk-j+1)}} \right)
+ \alpha_{2,d}m'.
\end{equation}

\noindent {\bf Incidences on the zero set $Z(f)$.} It remains to bound
incidences with points that lie on $Z(f)$. Set $\pts_0= \pts \cap Z(f)$
and $m_0=|\pts_0|=m-m'$. Let $\C_0$ denote the set of curves that are
fully contained in $Z(f)$, and set $\C' = \C \setminus \C_0$, $n_0 = |\C_0|$,
and $n' = |\C'| =n-n_0$. Since every curve of $\C'$ intersects $Z(f)$ in
$O(r^{1/d})$ points, we have, arguing as above,
\begin{equation} \label{eq:inc'd}
I(\pts_0,\C') \le b_dn'r^{1/d} = O(nr^{1/d}) \le \frac {\alpha_{1,d}}{3}
m^{\frac{k}{dk-d+1}+\eps}n^{\frac{dk-d}{dk-d+1}} ,
\end{equation}
provided that $\alpha_{1,d}$ is chosen sufficiently large.

Finally, we consider the number of incidences between points of
$\pts_0$ and curves of $\C_0$. For this, we set $c_{d-1}(k,d,\eps)$
to be the degree of $f$, which is $O(r^{1/d})=O((6b')^{1/(\eps
d)})$. Then, by the assumption of the theorem, we have $|\C_0| \le q_{d-1}$.
We consider a generic hyperplane $H \subset \RR^d$ and
project $\pts_0$ and $\C_0$ onto two respective sets $\pts^*$ and
$\C^*$ on $H$. Arguing as in the three-dimensional case, we can enforce
that $I(\pts_0,\C_0) = I(\pts^*,\C^*)$, that the projected curves have $k$ degrees of freedom,
%
%
%
and that, for $j<d-1$, the pairs $(q_j,c_j)$ remain
unchanged for $\pts^*$ and $\C^*$ within $H$.
Applying the induction hypothesis for dimension $d-1$, and recalling that
$|\C_0|\le q_{d-1}$, we obtain
\begin{equation*}
I(\pts_0,\C_0) = I(\pts^*,\C^*) \le \alpha_{1,d-1}
\left(\sum_{j=2}^{d-1}
m^{\frac{k}{jk-j+1}+\eps}q_{d-1}^{\frac{(d-1)(j-1)(k-1)}{(d-2)(jk-j+1)}}q_j^{\frac{(d-j-1)(k-1)}{(d-2)(jk-j+1)}}\right)
+\alpha_{2,d-1}(m+n).
\end{equation*}


As is easily verified, Equation~\eqref{eq:ontheq} with $l=d$ (and
$q_d = n$) implies that, for each $j$,
$$
q_{d-1}^{\frac{(d-1)(j-1)(k-1)}{(d-2)(jk-j+1)}}q_j^{\frac{(d-j-1)(k-1)}{(d-2)(jk-j+1)}}
\le n^{\frac{d(j-1)(k-1)}{(d-1)(jk-j+1)}}q_j^{\frac{(d-j)(k-1)}{(d-1)(jk-j+1)}}.
$$
By choosing $\alpha_{1,d}\ge 3\alpha_{1,d-1}$ and $\alpha_{2,d}\ge \alpha_{2,d-1}$, we have
that $I(\pts_0,\C_0)$ is at most
\begin{equation}
\frac{\alpha_{1,d}}{3}
\left(\sum_{j=2}^{d-1}
m^{\frac{k}{jk-j+1}+\eps}n^{\frac{d(j-1)(k-1)}{(d-1)(jk-j+1)}}q_j^{\frac{(d-j)(k-1)}{(d-1)(jk-j+1)}}\right)
+\alpha_{2,d}(m+n) . \label{eq:inc0d}
\end{equation}
By combining \eqref{eq:incCellsd}, \eqref{eq:inc'd}, and
\eqref{eq:inc0d}, including the case $m=O(n^{1/k})$,
and choosing $\alpha_{2,d}$ sufficiently large, we obtain
\begin{equation*}
I(\pts,\C) \le \alpha_{1,d} \left( m^{\frac{k}{dk-d+1} +\eps}n^{\frac{dk-d}{dk-d+1}}
+ m^{\frac{k}{2k-1}+\eps}n^{\frac{dk-d}{(d-1)(2k-1)}}q^{\frac{(k-1)(d-2)}{(d-1)(2k-1)}}\right)
+ \alpha_{2,d}(m+n).
\end{equation*}
This completes the induction step and thus the proof of the theorem.
\proofend

As a consequence of Theorem~\ref{th:IncD}, we have:

\noindent{\bf Example: incidences between points and lines in $\reals^4$.}

In the earlier version~\cite{SS14} of our study of point-line incidences in four dimensions,
we have obtained the following weaker version of Theorem~\ref{th:ss4d}.
\begin{theorem} \label{th:socg14}
For each $\eps>0$, there exists an integer $c_\eps$, so that the
following holds. Let $P$ be a set of $m$ distinct points and $L$ a
set of $n$ distinct lines in $\reals^4$, and let $q, s\le n$ be
parameters, such that
(i) for any polynomial $f \in \reals[x,y,z,w]$ of degree $\le
c_\eps$, its zero set $Z(f)$ does not contain more than $q$ lines of
$L$, and (ii) no 2-plane contains more than $s$ lines of $L$. Then,
$$
\label {ma:in}
I(P,L) \le A_{\eps} \left(m^{2/5+\eps}n^{4/5}+ m^{1/2+\eps}n^{2/3}q^{1/12}
+ m^{2/3+\eps}n^{4/9}s^{2/9} \right) + A(m + n) ,
$$
where $A_{\eps}$ depends on $\eps$, and $A$ is an absolute constant.
\end{theorem}

This result follows from our main Theorem~\ref{th:IncD}, if we impose
Equation~\eqref{eq:ontheq} on $q_2=s$, $q_3=q$, and $n$, which in this case is
equivalent to $s\le q \le n$ and $\frac {q^9} {n^8} <s$.
This illustrates how the general theory developed in this paper extends
similar results obtained earlier for ``isolated'' instances. Nevertheless,
as already mentioned earlier, the bound for lines in $\RR^4$ has been
improved in Theorem~\ref{th:ss4d} of \cite{SS4d}, in its
lower-dimensional terms.

\noindent{\bf Discussion.} We first notice that similarly to the
three-dimensional case, Theorem~\ref{th:IncD} implies an upper bound
on the number of $k$-rich points in $d$ dimensions (see
Corollary~\ref{co:rich} in three dimensions), and the proof thereof
applies verbatim, with the appropriate modifications of the various
exponents that now depend also on $d$. We leave it to the reader to
work out the precise statement.

Second, we note that Theorems \ref{th:mainInc3} and \ref{th:IncD}
have several weaknesses. The obvious ones are the items (i)--(iii)
discussed in the introduction. Another, less obvious weakness, has
to do with the way in which the $q_j$-dependent terms in the bounds
are derived. Specifically, these terms facilitate the induction
step, when the constraining parameters $q_j$ are passed
\emph{unchanged} to the inductive subproblems. Informally, since the
overall number of lines in a subproblem goes down, one would expect
the various parameters $q_j$ to decrease too, but so far we do not
have a clean mechanism for doing so. This weakness is manifested,
e.g., in Corollary \ref{co:rich}, where one would like to replace
the second term by one with a smaller exponent of $n$ and a larger
one of $q=q_2$. Specifically, for lines in $\RR^3$, one would like
to get a term close to $O(nq_2/k^3)$. This would yield
$O(n^{3/2}/k^3)$ for the important special case $q_2=O(n^{1/2})$
considered in \cite{GK2}; the present bound is weaker.

A final remark concerns the relationships between the parameters
$q_j$, as set forth in Equation~\eqref{eq:ontheq}. These conditions
are forced upon us by the induction process. As noted above, for
incidences between points and lines in $\reals^4$, the bound derived
in our main theorem~\ref{th:IncD} is (asymptotically) the same as
that of the main result of Sharir and Solomon in~\cite{SS14}. The
difference is that there, no restrictions on the $q_j$ are imposed.
The proof in~\cite{SS14} is facilitated by the so called ``second
partitioning polynomial'' (see ~\cite{KMSS,SS14}). Recently, Basu
and Sombra~\cite{BS14} proved the existence of a third partitioning
polynomial (see~\cite[Theorem 3.1]{BS14}), and conjectured the
existence of a $k$-partitioning polynomial for general $k>3$
(see~\cite[Conjecture 3.4]{BS14} for an exact formulation); for
completeness we refer also to ~\cite[Theorem 4.1]{FPSSZ14}, where a
weaker version of this conjecture is proved. Building upon the work
of Basu and Sombra~\cite{BS14}, the proof of Sharir and
Solomon~\cite{SS4d} is likely to extend and yield the same bound as
in our main Theorem~\ref{th:IncD}, for the more general case of
incidences between points and bounded degree algebraic curves in
dimensions at most five, and, if ~\cite[Conjecture 3.4]{BS14} holds,
in every dimension, without any conditions on the $q_j$.


\end{document}